 


\magnification=1200
\nopagenumbers
\parindent= 15pt
\baselineskip=13pt

\hsize=13.65cm
\vsize=19.3cm
\hoffset=-0.15cm
\voffset=0.3cm

\input amssym.def
\input amssym.tex

 at6.5pt
\font\srm=cmr8 

\font\cmx=cmbxti10
\font\csc=cmcsc10
\font\title=cmr12 at 16pt

\font\teneusm=eusm10    
\font\seveneusm=eusm7  
\font\fiveeusm=eusm5    
\newfam\eusmfam
\def\eusm{\fam\eusmfam\teneusm}
\textfont\eusmfam=\teneusm 
\scriptfont\eusmfam=\seveneusm
\scriptscriptfont\eusmfam=\fiveeusm

\def\Re{{\rm Re}\,}
\def\Im{{\rm Im}\,}
\def\sgn{{\rm sgn}\,}
\def\txt#1{{\textstyle{#1}}}
\def\scr#1{{\scriptstyle{#1}}}

\def\r#1{{\rm #1}}
\def\e#1{{\eusm #1}}
\def\varGamma{{\mit\Gamma}}
\def\B#1{{\Bbb #1}}
\def\b#1{{\bf #1}}
\def\sgn{{\rm sgn}}

\def\rightheadline{\hfil{\srm Extension of
the Linnik Phenomenon. II}
\hfil\tenrm\folio}
\def\leftheadline{\tenrm\folio\hfil{\srm Y. Motohashi}\hfil}
\def\emptyheadline{}
\headline{\ifnum\pageno=1 \emptyheadline\else
\ifodd\pageno \rightheadline \else \leftheadline\fi\fi}

\def\firstpage{\hss{\vbox to 1.0cm
{\vfil\hbox{\rm\folio}}}\hss}
\def\emptyfootline{\hfil}
\footline{\ifnum\pageno=1\firstpage\else
\emptyfootline\fi}

\centerline{\title An Extension of the
Linnik Phenomenon. II} 
\vskip 0.7cm
\centerline{\csc Yoichi Motohashi}
\vskip 1cm 
\noindent
{\bf Abstract:}  This work is a continuation of [9] but can be
read independently.
We discuss the extension of the Linnik phenomenon 
to automorphic $L$-functions
associated with cusp forms, focusing our attention on
the real analytic situation, as the holomorphic case
is settled in [9]. Our main assertion, which
is given at the end of the fifth section, reveals that
the repelling effect of exceptional zeros of
Dirichlet $L$-functions should be felt not only by
those $L$-functions themselves but also by automorphic
$L$-functions. We stress that constants, 
including those implicit, are all universal and
effectively computable, unless
otherwise stated. 
\bigskip
\noindent
{\bf 1. Introduction.}  
\par
\noindent
We repeat first the notion 
of exceptional zeros which we adopted in [9]: 
Let $\chi$ denote 
a generic Dirichlet character, 
with which the $L$-function 
$L(s,\chi)$ is associated. Let $Z_T=\{\rho=
\beta+i\gamma\}$ be the set of all non-trivial zeros $\rho$
in the region $|\Im s|\le T$ of the function
$\prod_{q\le T}\prod^*_{\chi\bmod q}L(s,\chi)$; 
here and in what follows $T>0$ 
is assumed to be sufficiently large and
the asterisk means 
that relevant characters are primitive. Then, there exists a
constant $a_0>0$ such that
$$
\hbox{$\displaystyle
\max_{\rho\in Z_T}\beta< 1-{a_0\over\log T},\;$
save for a possibly existing $\beta_T\in Z_T$.
}\eqno(1.1)
$$
If $\beta_T$ ever exists, it should be real and simple, 
and we designate both itself and
the relevant unique primitive character $\chi_T$, with
$L(\beta_T,\chi_T)=0$, 
as $T$-exceptional. It is known also that $\chi_T$ 
should be real. As a matter of fact, we have more precisely
$$
1-{a_0\over \log T}\le\beta_T
\le 1-{1\over \sqrt{T}(\log T)^3};\eqno(1.2)
$$
for a proof see [7, Vol.\ I, Chapter 4] for instance.
\par
Superseding the upper bound in $(1.2)$, 
a theorem of Siegel asserts that
for any fixed $\varepsilon>0$ there exists a 
$c(\varepsilon)>0$ such that
$$
\beta_T\le1-{c(\varepsilon)\over T^\varepsilon}.
\eqno(1.3)
$$
However, all known proofs of $(1.3)$ yield only 
the existence of $c(\varepsilon)$ and do not provide 
any means to evaluate its actual values. 
This ineffectiveness in $(1.3)$ causes severe difficulties 
in various basic problems. The most outstanding
among them is the estimation of  the size 
of the least prime that appears in a 
given  arithmetic progression. 
To resolve this particular difficulty,  Linnik
greatly refined $(1.3)$ by providing
a  quantitative version of 
the Deuring--Heilbronn phenomenon or
the repelling effect of the existence of $\beta_T$
towards all other elements in $Z_T$. Linnik's
theorem  or rather his phenomenon asserts that the inequality
$$
\beta\le1-{a_0\over\log T}\log{a_0 e\over
(1-\beta_T)\log T}\eqno(1.4)
$$
holds for all $\rho\ne\beta_T$ in $Z_T$.
With this and a zero density result of a special type 
for Dirichlet $L$-functions,
which is another basic contribution of his, 
Linnik could reach his famed Least Prime Number 
Theorem for arithmetic progressions. A brief history of the
research relevant to the theorem
is given in [7, Vol.\ I, Chapter 9] and [9, Section 1]. 
\medskip
In the present work we shall show a way to
extend (1.4) to zeros of 
automorphic $L$-functions associated with
real analytic cusp forms; namely, $(1.4)$ holds even if we
include into $Z_T$ all the zeros of automorphic $L$-functions
whose analytic conductors are of polynomial order in $T$.
It should be stressed that the holomorphic case
is successfully resolved in [9].
The difference between the holomorphic and the real analytic
situations in our context is in that 
the Ramanujan conjecture is confirmed in the former
case but not yet in the latter. Thus 
we have to devise a new approach in dealing
with the real analytic case.  An important key to overcome
this difficulty lies in a structure of 
Rankin $L$-functions discovered 
by Shimura [12] and extended by Gelbart--Jacquet [1]
in a far-reaching manner.
\medskip
We shall restrict ourselves to the framework 
$$
\r{G}=\r{PSL}(2,\B{R}),\quad
\varGamma=\r{PSL}(2,\B{Z}).\eqno(1.5)
$$
This is solely for the sake of simplicity. Extension to
congruence groups should not cause any extra difficulty.
A fairly elementary theory of automorphic 
representations within $(1.5)$
is developed in [7, Vol.\ II, Chapters\ 2--3] and [8].
\medskip
\noindent
{\csc Acknowledgements}: 
The author is indebted to S. Gelbart, 
J. Hoffstein, A. Perelli, A. Sankaranarayanan
and N. Watt for their kind assists and comments.
\bigskip
\noindent
{\bf 2. Automorphic forms.}
\par
\noindent
We need to make precise how we normalise 
$\varGamma$-automorphic forms.
We begin with some basic notions on $\r{G}$. 
Each $\r{h}\in\r{G}$ induces the left
and the right translations $l_\r{h}(\r{g})=\r{hg}$,
$ r_\r{h}(\r{g})=\r{gh}$, which are real analytic
with respect to the Iwasawa co-ordinate system
$\r{G}=\r{NAK}\ni\r{g}=\r{n}[x]\r{a}[y]\r{k}[\theta]$,
with
$$
\displaystyle{\r{N}=\left\{\r{n}[x]=
\left[\matrix{1&x\cr&1}\right]:
\,x\in{\Bbb R}\right\},\quad\r{A}=\left\{\r{a}[y]
=\left[\matrix{\sqrt{y}&\cr&1/\sqrt{y}}\right]:
\,y>0\right\},}\atop\displaystyle{
\r{K}=\left\{\r{k}[\theta]=
\left[\matrix{\phantom{-}\cos\theta&
\sin\theta\cr-\sin\theta&\cos\theta}\right]:\,\theta\in{\Bbb
R}/\pi{\Bbb Z}\right\}.
}\eqno(2.1)
$$
The symbols $\{x,y,\theta\}$ will retain this specification
within the present section.
We put $d\r{g}={dxdyd\theta/\pi y^2}$
with Lebesgue measures $\, dx, dy, d\theta$. 
The group $\r{G}$ is unimodular, as
we have
$d\r{g}=dl_\r{h}(\r{g})$, 
$d\r{g}=dr_\r{h}(\r{g})$, $\forall\r{h}\in\r{G}$.
The centre of all left invariant differential operators 
on $\r{G}$ is the polynomial ring generated 
by the Casimir operator
$$
\Omega=-y^2\left(\Big({\partial\over\partial x}\Big)^2
+\Big({\partial\over\partial y}\Big)^2\right)
+y{\partial^2\over\partial x\partial\theta}.\eqno(2.2)
$$
This is not only left but also right invariant.
\par
Next, let $L^2(\varGamma\backslash\r{G})$ be the set of all
functions $f$ or vectors on $\r{G}$ which are left
$\varGamma$-automorphic, i.e., $l_\gamma f=f$, 
$\forall\gamma\in\varGamma$,
and square integrable against $d\r{g}$ 
over any fundamental 
domain of $\varGamma$ on $\r{G}$;
for instance, the set $[\e{F}]=\big\{\r{n}[x]
\r{a}[y]\r{k}[\theta]:
x+iy\in\e{F},\, 0\le\theta\le\pi\big\}$,
with $\e{F}=\{z\in\B{C}:|z|\ge 1, 
|\Re z|\le{1\over2}\}$, is such a domain. Then
$L^2(\varGamma\backslash\r{G})$ is 
a Hilbert space equipped with the Petersson inner-product
$$
\big\langle f_1, f_2\big\rangle
=\int_{\varGamma\backslash\r{G}}
f_1(\r{g})\overline{f_2(\r{g})}d\r{g},\eqno(2.3)
$$
where the integration range is the whole quotient space 
$\varGamma\backslash\r{G}$ and the measure is the
one naturally induced. Then, let $\Vert{f}\Vert$ 
be the norm associated with $(2.3)$. The unimodularity 
of $d\r{g}$ implies that right translations are all unitary 
maps of $L^2(\varGamma\backslash\r{G})$ onto
itself. The map $r:\r{h}\mapsto r_\r{h}$, which
is a strongly continuous homomorphism of $\r{G}$ into the
unitary transformation group of 
$L^2(\varGamma\backslash\r{G})$, is termed the regular 
$\varGamma$-automorphic representation of $\r{G}$. 
Any closed subspace $W$ of 
$L^2(\varGamma\backslash\r{G})$ which satisfies
$r_\r{h}W\subseteq W$, $\forall\r{h}\in\r{G}$, is 
called an invariant subspace and induces a unitary
representation of $\r{G}$. In this context we use
a representation and an invariant subspace as interchangeable
notions.
If $W$ does not contain any non-trivial invariant subspace,
then it is said to be an irreducible subspace or 
representation.
\par
In terms of the regular 
automorphic representation the space
$L^2(\varGamma\backslash\r{G})$ is spectrally decomposed
into irreducible subspaces. 
We have first the orthogonal decomposition 
into invariant subspaces
$$
L^2(\varGamma\backslash\r{G})={\Bbb C}\cdot1\oplus
{}^c\!L^2(\varGamma\backslash\r{G})
\oplus{}^e\!L^2(\varGamma\backslash\r{G}).\eqno(2.4)
$$
Here
${}^c\!L^2(\varGamma\backslash\r{G})$ is the 
cuspidal subspace which is
spanned by all vectors whose constant terms in the
Fourier expansion with respect to the left action of
$\r{N}$ vanish. 
The subspace ${}^e\!L^2(\varGamma\backslash\r{G})$ 
is generated by integrals of Eisenstein series; this part is
irrelevant to our present purpose. 
The spectral structure of the cuspidal
subspace is embodied in
$$
{}^c\!L^2(\varGamma\backslash\r{G})
=\oplus V,\eqno(2.5)
$$
where $V$ are all irreducible. At each $V$ the Casimir 
operator becomes a
constant multiplication or $V$ is an eigen-space 
of $\Omega$, which we express as
$$
\Omega|_V=\left(\txt{1\over4}-\nu_V^2\right)
\cdot 1,\quad\nu_V\in\B{C}.\eqno(2.6)
$$
In fact it should be
made precise in which domain we consider the action of
$\Omega$, but we skip the details. 
Because of the unimodularity of $d\r{g}$
the Casimir operator is symmetric; and 
${1\over 4}-\nu_V^2$ should be in $\B{R}$.
\par
There are two main series of irreducible cuspidal 
representations:
$$
\eqalign{
&\hbox{$V$ in the unitary principal series $\Leftrightarrow$
$\nu_V\in i\B{R}$,}\cr
&\hbox{$V$ in the discrete series $\Leftrightarrow$
$\nu_V\in \B{N}-{1\over2}$.}
}\eqno(2.7)
$$
In general there is also the complementary series where
$0\le\nu_V<{1\over2}$; but under the supposition 
$\varGamma=\r{PSL}(2,\B{Z})$ such does not occur.
Our discussion will be focused on $V$'s 
in the unitary principal 
series, since they are in fact
generated by Maass waves or
real analytic cusp forms on the hyperbolic
upper half plane $\r{G}/\r{K}$, as is explained below.
The $V$'s in the discrete series 
are generated by holomorphic and anti-holomorphic cusp forms.
\par
We consider further the Fourier expansion of elements in any
$V$ belonging to the unitary principal series,
with respect to the right action of $\r{K}$. We have
$$
V=\mathop\oplus_{\ell=-\infty}^\infty V_\ell,\eqno(2.8)
$$
where $f\in V_\ell$ means that $f(\r{g}\r{k}[\tau])=
\exp(2i\ell\tau)f(\r{g})$.  
We term $f$ a 
$\varGamma$-automorphic cusp form of weight $2\ell$ 
with the spectral data $\nu_V\in i\B{R}$ or 
just an automorphic form.
With the operators $\b{e}^\pm=e^{\pm2i\theta}
\big(\pm2iy\partial_x+2y\partial_y\mp i\partial_\theta\big)$,
we have
that the weights of $\b{e}^\pm f$ equal $2(\ell\pm1)$. 
This mechanism is the Maass shift. 
To make the situation more explicit, we
introduce the Jacquet operator: For a function $\phi$ on
$\r{G}$, we put
$$
\e{A}^\epsilon\!\phi(\r{g})=\int_{-\infty}^\infty 
e(-\epsilon\xi)\phi(\r{w}\r{n}[\xi]\r{g})d\xi,\quad
\epsilon=\pm,\> \r{w}=\r{k}
\big[\txt{1\over2}\pi\big],\> e(\xi)=\exp(2\pi i\xi).
\eqno(2.9)
$$
Then we have
$$
V_\ell=\B{C}\cdot\psi_V^{(\ell)},\quad
\langle\psi_V^{(\ell_1)},\psi_V^{(\ell_2)}\rangle=
\delta_{\ell_1,\,\ell_2},\eqno(2.10)
$$
with the Kronecker delta and
$$
\psi_V^{(\ell)}(\r{g})=\sum_{\scr{n=-\infty}
\atop\scr{n\ne0}}
^\infty{\varrho_V(n)\over\sqrt{|n|}}
\e{A}^{\sgn(n)}\phi_\ell
(\r{a}[|n|]\r{g},\nu_V),\quad
\phi_\ell(\r{g},\nu)=y^{\nu+1/2}
\exp(2\ell i\theta). \eqno(2.11)
$$
What is essential in 
our normalisation $(2.10)$--$(2.11)$ of automorphic forms
is that 
the sequence $\{\varrho_V(n): \B{Z}\ni n\ne0\}$ 
does not depend on the weights but is 
uniquely determined by $V$ save for constant multipliers
of unit absolute value. 
\par
In literature the transform $\e{A}^\delta\phi_\ell(\r{g},\nu)$
is usually expressed in terms of the
Whittaker function $W_{\mu,\nu}$:
$$
\e{A}^\epsilon\phi_\ell(\r{g},\nu)
=(-1)^\ell\pi^{\nu+1/2}\exp(2\ell i\theta)e(\epsilon x)
{W_{\epsilon \ell,\nu}(4\pi y)\over\Gamma
\big(\nu+\epsilon \ell+{1\over2}\big)}.\eqno(2.12)
$$
Or more explicitly, we have
$$
\eqalignno{
\e{A}^\epsilon\phi_\ell(\r{g},\nu)
&={\pi^\nu\exp(2\ell i\theta)e(\epsilon x)
\over\Gamma\big(\nu+|\ell|+{1\over2}\big)}
\sum_{j=0}^{|\ell|}(-\pi)^j{2|\ell|\choose 2j}
\Gamma\big(|\ell|-j+\txt{1\over2}\big)\cr
&\times\, y^{j+1/2}\int_0^\infty u^{\nu-1}
\big(\sqrt{u}+\epsilon\sgn(\ell)/\sqrt{u}\big)^{2j}
\exp\left(-\pi y\big(u+1/u\big)\right)du.\qquad&(2.13)
}
$$
In particular, we have
$$
\psi_V^{(0)}(\r{g})={2\pi^{\nu_V+1/2}\over
\Gamma\big(\nu_V+{1\over2}\big)}
\sum_{n\ne0}\varrho_V(n)y^{1/2}K_{\nu_V}(2\pi|n|y)
e(nx),\eqno(2.14)
$$
with $K_\nu$ being the $K$-Bessel function of order $\nu$;
thus $\psi_V^{(0)}$ is in fact a real analytic cusp form 
on the hyperbolic upper half plane $\r{G}/\r{K}$.
The Maass shift is the same as the identity
$$
\psi_V^{(\ell)}={\Gamma({1\over2}+\nu)\over
2^{|\ell|}\Gamma\big({1\over2}+\nu+|\ell|\big)}
(\b{e}^{\sgn(\ell)})^{|\ell|}\psi_V^{(0)},
\quad \ell\in\B{Z},\eqno(2.15)
$$
which is a consequence of the obvious intertwining
property of $\e{A}^\epsilon$. In this context,
any irreducible subspace belonging to
the unitary principal series is generated by
a Maass wave.
\par
We may in fact deal with non-negative weights only,
because of the existence of the involution 
$J:\r{n}[x]\r{a}[y]\r{k}[\theta]
\mapsto \r{n}[-x]\r{a}[y]\r{k}[-\theta]$
which is the same as the map
$\left[{a\atop c}{b\atop d}\right]\mapsto
\left[{\hfill a\atop -c}{-b\atop\hfill d}\right]$ in $\r{G}$. 
Since $J$ commutes with left translations,
we may assume in particular that
$J\psi_V^{(0)}=\epsilon_V\psi_V^{(0)}$ 
with $\epsilon_V=\pm1$,
which is equivalent to 
$$
\varrho_V(-n)=\epsilon_V\varrho_V(n),
\quad n\in\B{N}.\eqno(2.16)
$$ 
\bigskip
\noindent
{\bf 3. Automorphic {\cmx L}-functions.}
\par
\noindent
We turn to $L$-functions 
associated with irreducible subspaces $V$
in the unitary principal series. Before introducing
them, we impose that $V$ be all
Hecke invariant as well. Thus, let
$$
T(n)f(\r{g})={1\over\sqrt{n}}\sum_{ad=n}\,
\sum_{b\bmod d}f\Big(\Big[\matrix{a&b\cr&d}\Big]
\r{g}\Big),\quad n\in\B{N},\eqno(3.1)
$$
and we may assume 
that each $V$
is an eigen-space of every $T(n)$ so that there exists 
a real sequence $\{\tau_V(n): n\in\B{N}\}$ such that
$T(n)|_V=\tau_V(n)\cdot 1$. In particular we have,
together with $(2.16)$,
$$
\varrho_V(n)=\epsilon_V^{(1-\sgn(n))/2}\varrho_V(1)
\tau_V(|n|),\quad 
n\in\B{Z}\backslash\{0\}.\eqno(3.2)
$$
According to Kim--Sarnak [4, Appendix],
we have the bound
$$
|\tau_V(n)|\le d(n)n^{7/64},\quad n\in\B{N},
\eqno(3.3)
$$
uniformly in $V$, where $d(n)$ is the number of divisors of
$n$.
\par
Then we define the Hecke $L$-function associated 
with $V$ by
$$
L(s;V)=\sum_{n=1}^\infty \tau_V(n)n^{-s}.\eqno(3.4)
$$
The uniform and absolute convergence for $\Re s>1$
follows from that of $(3.9)$ below. This admits
the Euler product expansion
$$
\eqalignno{
L(s;V)&\,=\prod_p\Big(1-{\tau_V(p)\over p^s}
+{1\over p^{2s}}\Big)^{-1}\cr
&\,=\prod_p\Big(1-{\alpha_V(p)\over p^s}\Big)^{-1}
\Big(1-{\alpha^{-1}_V(p)\over p^s}\Big)^{-1},&(3.5)
}
$$
where $1\le|\alpha_V(p)|\le p^{7/64}\,$; 
here and in what follows
$p$ stands for a generic prime.
We shall need analytic properties of
the twist of $H(s;V)$ by a primitive Dirichlet
character $\chi$ with the conductor $q\ge1$:  
$$
L(s,\chi;V)=\sum_{n=1}^\infty \chi(n)\tau_V(n)n^{-s}
\eqno(3.6)
$$
This continues to an entire function and satisfies
the functional equation of the Riemannian type:
$$
{\Gamma_2(s)L(s,\chi;V)=C^{(1)}_\chi
\Gamma_2(1-s)
L(1-s,\overline\chi;V),}\atop
\displaystyle{\Gamma_2(s)=(q/\pi)^{s}
\prod_{j=1}^2\Gamma\big(
\txt{1\over2}(s+c_j)\big),}\eqno(3.7)
$$
with $|C^{(1)}_\chi|=1$ 
and $\sum_{j=1}^2|c_j|\ll|\nu_V|$;
our notation for the $\Gamma$-factor is 
somewhat over-simplified one.
Thus the Phragm\'en--Lindel\"of 
convexity principle and $(3.3)$ yield the uniform bound
$$
L(\sigma+it,\chi;V)\ll (q(|t|+|\nu_V|))^{10/9-\sigma},
\quad -\txt{1\over9}\le\sigma\le\txt{10\over 9}.
\eqno(3.8)
$$
A better bound is of course available, but this is more than
enough for our purpose.
\par
Also we shall require basic properties 
of the Rankin $L$-functions associated with $V$
and their twists:
$$
\eqalign{
L(s;V\times V)&=\zeta(2s)
\sum_{n=1}^\infty\tau_V^2(n)n^{-s},\cr
L(s,\chi;V\times V)&=L(2s,\chi^2)
\sum_{n=1}^\infty\chi(n)
\tau_V^2(n)n^{-s},}\eqno(3.9)
$$
with $\chi$ as above.
Via the unfolding method, 
the absolute convergence for $\Re s>1$ of the first follows. 
In fact, $L(s; V\times V)$ is regular except for the
simple pole at $s=1$ with the residue 
${1\over2}|\varrho_V(1)|^{-2}$, which is but a fact 
that we shall not use. The function
$L(s,\chi;V\times V)$ is 
entire for $q>1$ and that for $q\ge1$
the functional equation
$$
{\Gamma_4(s)L(s,\chi;V\times V)=C^{(2)}_\chi
\Gamma_4(1-s)R(1-s,\overline{\chi};V\times V),}\atop
\displaystyle{\Gamma_4(s)=(q/\pi)^{2s}
\prod_{j=1}^4\Gamma\big(\txt{1\over2}
(s+d_j)\big)}\eqno(3.10)
$$
holds; here $|C^{(2)}_\chi|=1$ and 
$\sum_{j=1}^4 |d_j|\ll|\nu_V|$; the same simplification
as in $(3.7)$ is adopted.
In particular,
the convexity argument gives
$$
L(\sigma+it,\chi;V\times V)\ll \big(q(|t|
+|\nu_V|)\big)^{2(1+\eta-
\sigma)}|\nu_V|^{4/9-2\eta},
\quad -\eta\le\sigma\le 1+\eta,
\eqno(3.11)
$$
for any fixed $\eta\in\big(0,{2\over9}\big]$, with
the understanding that
if $q=1$ then a small neighbourhood of $s=1$ is to be 
excluded. In fact, $(3.3)$ gives $L(\sigma;V\times V)
\ll1$ for $\sigma\ge{11\over9}$, which yields 
$|L(1+\eta+it,\chi;V\times V)|\le L(1+\eta;V\times V)
\ll |\nu_V|^{4/9-2\eta}$ via the
absolute convergence and the convexity. Then, applying
the convexity argument again, we get $(3.11)$. A better
bound is available but $(3.11)$ serves well for our purpose.
\par
We remark also that $(3.5)$ implies
the following relations: In the region of absolute convergence,
$$
\eqalignno{
\sum_{\ell=0}^\infty{\sigma_a(p^\ell,\chi)\tau_V(p^\ell)
\over p^{\ell s}}&\,=\Big(1-{\chi(p)\over p^{2s-a}}\Big)
\Big(1-{\alpha_V(p)\over p^s}\Big)^{-1}
\Big(1-{\alpha^{-1}_V(p)\over p^s}\Big)^{-1}\cr
&\,\times\Big(1-{\chi(p)\alpha_V(p)\over p^s}
\Big)^{-1}\Big(1-{\chi(p)\alpha^{-1}_V(p)\over p^s}
\Big)^{-1},&(3.12)
}
$$
with $\sigma_a(n,\chi)=\sum_{d|n}\chi(d)d^a$; and
$$
\eqalignno{
\sum_{\ell=0}^\infty{\chi(p^\ell)\tau_V^2(p^\ell)\over
p^{\ell s}}&\,=\Big(1-{\chi(p^2)\over p^{2s}}\Big)
\Big(1-{\chi(p)\over p^s}\Big)^{-2}\cr
&\,\times\Big(1-{\chi(p)\alpha_V^2(p)\over p^s}\Big)^{-1}
\Big(1-{\chi(p)\alpha_V^{-2}(p)\over
p^s}\Big)^{-1}.&(3.13)
}
$$
\medskip
The next assertion, which stems from
Gelbart--Jacquet [1, (9.3) Theorem], 
is the most basic structure of 
Rankin $L$-functions that we exploit:
\medskip
\noindent
{\bf Lemma 1.} {\it There exists an
entire $L$-functions $\e{L}(s;V)$ and 
$\e{L}(s,\chi;V)$ such that
$$
\eqalign{
L(s;V\times V)&\,=\zeta(s)\e{L}(s;V),\cr
L(s,\chi;V\times V)&\,=L(s,\chi)\e{L}(s,\chi; V).
}\eqno(3.14)
$$ 
\/}
\smallskip
\noindent 
Hence twists of Rankin $L$-functions may have
exceptional zeros, i.e., 
$L(\beta_T,\chi_T;V\times V)=0$. Our
argument is based on this particular fact. We should mention
also that $\e{L}(s;V)$ is usually
denoted as $L(s;\r{sym}^2V)$
the symmetric 2nd power $L$-function associated with
the representation $V$.
\bigskip
\noindent
{\bf 4. Sieve tools.}
\par
\noindent
We collect basic sieve devices which are to be employed 
in the next section. Hereafter we shall assume 
that the $T$-exceptional character $\chi_T$ exists. 
We put $\beta_T=1-\kappa$, and introduce 
the multiplicative function
$$
f(n)=\sigma_{-\kappa}(n,\chi_T)=
\sum_{d|n}\chi_T(d)d^{-\kappa},
\eqno(4.1)
$$
which is always positive. Also we are concerned
with those irreducible representations $V$ 
in the unitary principal series such that
$$
\log |\nu_V|\ll \log T.\eqno(4.2)
$$
With this, we consider the 
$\Lambda^2$-sieve situation
$$
\sum_{n\le N}f(n)\tau_V^2(n)
\Bigg(\sum_{d|n}\lambda_d\Bigg)^2,
\eqno(4.3)
$$
where the real numbers $\{\lambda_d\}$ are supported
on the set of square-free integers and such that
$\lambda_1=1$, and $\lambda_d=0$ for $d>R\ge1$. 
To set out
our choice of $\{\lambda_d\}$ we introduce
$$
F_p(s)=\sum_{\ell=0}^\infty f(p^\ell)\tau_V^2(p^\ell)
p^{-\ell s},\eqno(4.4)
$$
which converges absolutely for $\Re s>{7\over32}$ 
because of $(3.3)$. We put $F(s)=\prod_p
F_p(s)$.  We have
$$
F(s)=L(s;V\times V)
L(s+\kappa,\chi_T;V\times V)Y(s),\eqno(4.5)
$$
where $Y(s)$ is regular and bounded for 
$\Re s>{23\over 32}$. In fact,
by definition, we have, for $\Re s>1$,
$$
Y(s)=\prod_p{\big(1-p^{-2s}\big)\big(1-\chi_T^2(p)
p^{-2(s+\kappa)}\big)F_p(s)\over
\displaystyle\bigg(\sum_{\ell=0}^\infty\tau_V^2(p^\ell)
p^{-\ell s}\bigg)\bigg(\sum_{\ell=0}^\infty
\chi_T(p^\ell)\tau_V^2(p^\ell)
p^{-\ell(s+\kappa)}\bigg)}.\eqno(4.6)
$$
We divide this into two parts
 $Y=Y_1Y_2$ where $Y_1$ is the product over $p\le p_0$
and $Y_2$ the rest. With any finite
$p_0$, the part $Y_1$ is regular
for $\Re s> {7\over 32}$ 
because of $(3.3)$ and $(3.13)$. On the other hand, 
if $p_0$ is large, then the $p$-factor in $Y_2$ is $1+
O\big(\sum_{\ell=2}^\infty(\ell+1)^5
p^{(7/32-\Re s)\ell}\big)$. Hence
$Y_2$ is regular and bounded 
in the range indicated above.
\par
Then, with $F_p=F_p(1)$, we put
$$
\lambda_d=\mu(d)F_d{G_d(R/d)\over G_1(R)},
\eqno(4.7)
$$
where $\mu$ is the M\"obius function, and
$$
F_d=\prod_{p|d}F_p,\quad
G_d(x)=\sum_{\scr{r\le x}\atop\scr{(d,r)=1}}
{\mu^2(r)K(r)},\quad K(r)=\prod_{p|r}
(F_p-1),\eqno(4.8)
$$
with the greatest common divisor $(d,r)$ of $d$ and $r$.
It should be stressed that $K(r)>0$ always. 
In fact we have more explicitly that
$$
F_p-1\ge {1\over 16p^{4(1+\kappa)}}.\eqno(4.9)
$$
To show this we note that
$$
F_p-1\ge{\tau^2_V(p^2)\over p^2}f(p^2)
+{\tau^2_V(p^4)\over p^4}f(p^4)
\ge{\tau^2_V(p^2)\over p^{2(1+\kappa)}}
+{\tau^2_V(p^4)\over p^{4(1+\kappa)}},\eqno(4.10)
$$
and that
$$
\tau_V(p^4)=\big(\tau_V(p^2)-\txt{1\over2}\big)^2
-\txt{5\over4}.\eqno(4.11)
$$
If $|\tau_V(p^4)|\ge{1\over4}$, then $(4.9)$ is
immediate. Otherwise, $|\tau_V(p^2)-{1\over2}|\ge1$
and thus $|\tau_V(p^2)|\ge{1\over2}$, which
implies $(4.9)$.
\medskip
Now, the choice $(4.7)$ 
leads us to the multiplicative function
$$
\Phi_r(n)={\mu((n,r))\over K((n,r))};\eqno(4.12)
$$
see either [6, \S1.4] or [7, Vol.\ I, Chapter 9]. 
We are about to show the
quasi-orthogonality in the set $\{\Phi_r(n)\}$.
To this end we consider the expression
$$
\sum_{n\le N}f(n)\tau^2_V(n)
\Bigg|\sum_{r\le R}
\mu^2(r)\Phi_r(n)\sqrt{K(r)}\cdot
b_r\Bigg|^2,\eqno(4.13)
$$
where $N, R\ge1$ and $\{b_r\}$ are arbitrary.
Expanding the squares out, we have
$$
\sum_{r_1,\,r_2\le R}{\mu^2(r_1)\mu^2(r_2)
\sqrt{K(r_1)K(r_2)}}\,S(N;r_1,r_2)
b_{r_1}\overline{b}_{r_2},\eqno(4.14)
$$
with
$$
S(N; r_1,r_2)=
\sum_{n\le N}f(n)\tau_V^2(n)
\Phi_{r_1}(n)\Phi_{r_2}(n).\eqno(4.15)
$$
Thus let us consider the function
$$
\eqalignno{
&\sum_{n=1}^\infty f(n)\tau_V^2(n)
\Phi_{r_1}(n)\Phi_{r_2}(n)n^{-s}\cr
=&\,\Bigg(\sum_{(n,\, r_1r_2)=1}\Bigg)
\Bigg(\sum_{n|([r_1,\,r_2]/(r_1,\,r_2))^\infty}\Bigg)
\Bigg(\sum_{n|( r_1,\,r_2)^\infty}\Bigg)
=\r{F}_1\r{F}_2\r{F}_3,&(4.16)
}
$$
say, where $[r_1,r_2]$ is the least common multiple of
$r_1$ and $r_2$, and it is assumed that $\Re s$ is
sufficiently large. We have
$$
\eqalign{
\r{F}_1&=F(s)\prod_{p\mid r_1r_2}F_p(s)^{-1},\cr
\r{F}_2&=\prod_{p|([r_1,\,r_2]/(r_1,\,r_2))}
\big(1-(F_p-1)^{-1}(F_p(s)-1)\big),\cr
\r{F}_3&=\prod_{p|(r_1,\,r_2)}
\big(1+(F_p-1)^{-2}(F_p(s)-1)\big),
}\eqno(4.17)
$$
which gives
$$
{\r{F}_1\r{F}_2\r{F}_3=F(s)U_{r_1,r_2}(s),}\atop
\displaystyle{U_{r_1,r_2}(s)=\r{F}_2\r{F}_3
\prod_{p|r_1r_2}F_p(s)^{-1}
=\sum_{d|(r_1r_2)^\infty} u(d)d^{-s},
}\eqno(4.18)
$$
say; that is,
$$
f(n)\tau_V^2(n)\Phi_{r_1}(n)
\Phi_{r_2}(n)=\sum_{d|(n,(r_1r_2)^\infty)}f(n/d)
\tau_V^2(n/d)u(d).\eqno(4.19)
$$
Thus
$$
S(N;r_1,r_2)=\sum_{d|(r_1r_2)^\infty} u(d)
\sum_{n\le N/d}f(n)\tau_V^2(n),\eqno(4.20)
$$
where empty sums are to vanish.
To evaluate the inner-sum we use
Perron's inversion formula. We multiply $(4.5)$ by
$N^s/(2\pi is)$ and integrate over the segment
$[1+\eta-iP, 1+\eta+iP]$, where $\eta$ is as in
$(3.11)$ and $P\ge1$ to be fixed later. We move the
segment to $\Re s={3\over 4}$, encountering only one
singularity at $s=1$ because of Lemma 1. We get,
on $(4.2)$,
$$
\eqalignno{
\sum_{n\le N}f(n)\tau_V^2(n)
=&\,\e{L}(1;V)L(1+\kappa,\chi_T;V\times V)Y(1)N\cr
&\qquad+\,O\big(T^c(N^{39/32+\eta}/P
+N^{3/4}P^{1+4\eta})\big)\cr
=&\,\e{L}(1;V)L(1+\kappa,\chi_T;V\times V)Y(1)N
+O\big(T^c N^{64/65}\big),&(4.21)
}
$$
with a constant $c>0$. Inserting this into $(4.20)$, we 
find that on $(4.2)$
$$
\eqalignno{
S(N;r_1,r_2)=&\,
{\delta_{r_1,r_2}\over K(r_1)}\e{L}(1;V)
L(1+\kappa,\chi_T;V\times V)Y(1)N\cr
&+O\big(T^c (r_1r_2)^{9/2}N^{64/65}\big),&(4.22)
}
$$
uniformly for $N\ge1$ and
square-free $r_1, r_2\ge1$,
since a combination of $(4.9)$ and $(4.17)$--$(4.18)$
gives
$$
\sum_{d|(r_1r_2)^\infty}|u(d)|d^{-64/65}\le
D^{\nu(r_1r_2)}(r_1r_2)^{4(1+\kappa)},\eqno(4.23)
$$
where $D$ is a constant and $\nu(r)=\sum_{p|r}1$;
in fact one may show a much better bound but this suffices
for our purpose. Collecting these assertions and
invoking the duality principle, we obtain the
following analogue of [6, Theorem 5][7, Vol.\ I, (9.1.23)]:
\medskip
\noindent
{\bf Lemma 2.} {\it 
There exists a positive constant $c$ such that
we have, for any $N,R\ge1$ and $\{a_n\}$,
$$
\eqalignno{
&\sum_{r\le R}\mu^2(r)K(r)
\Bigg|\sum_{n\le N}f(n)\tau_V(n)
\Phi_r(n)a_n\Bigg|^2\cr
\ll &\,\Big(\e{L}(1;V)L(1+\kappa,\chi_T;V\times V)N
+T^cR^{10}N^{64/65}\Big)
\sum_{n\le N}f(n)|a_n|^2,&(4.24)
}
$$
uniformly for $V$ satisfying $(4.2)$.
\/}
\medskip
\noindent
This exponent $64\over65$ is by no means the best
that one can attain; any explicit exponent 
less than $1$ works fine.
\medskip
On the other hand, the sieve effect of $(4.7)$ is embodied in 
\medskip
\noindent
{\bf Lemma 3.} {\it On $(4.2)$ we have
$$
G_1(R)\asymp\kappa^{-1}\e{L}(1;V)
L(1+\kappa,\chi_T;V\times V),\eqno(4.25)
$$
provided 
$$
\hbox{$\log R/\log T$ is sufficiently large but bounded.}
\eqno(4.26)
$$
In particular we have
$$
\e{L}(1;V)L(1+\kappa,\chi_T;V\times V)\gg 
\big(\sqrt{T}(\log T)^3\big)^{-1}.\eqno(4.27)
$$
}
\par
\noindent
{\it Proof\/}. The bound $(4.25)$ 
corresponds to [6, (6.2.5)][7, Vol.\ I, (9.2.51)]; 
$\kappa$ here
is the same as $\delta$ there. We note that
under $(4.26)$
$$
G_1(R)\asymp R^{-2\kappa}\sum_{r\le R}\mu^2(r)
K(r)r^{2\kappa},\eqno(4.28)
$$
and that
$$
\sum_{r=1}^\infty \mu^2(r)K(r)r^{-s+2\kappa}
=L(s+1-2\kappa;V\times V)
L(s+1-\kappa,\chi_T; V\times V)W(s),
\eqno(4.29)
$$
where $W(s)$ is regular and bounded for 
$\Re s>-{9\over32}$, which can be shown 
in much the same way as $(4.5)$.
We multiply this by
$R^s/(2\pi i s)$ and integrate over the vertical segment
$[\eta-iP,\eta+iP]$ with a small constant 
$\eta>2\kappa$ and 
with $P\ge1$ to be fixed later. 
Then we shift the segment to $\Re s=-{1\over4}$,
encountering only one pole at $s=2\kappa$, 
since the integrand is
regular at $s=0$ because of Lemma 1. Skipping
the remaining details, since they are the same 
as those for $(4.21)$,
we assert that
$$
\sum_{r\le R}\mu^2(r)K(r)r^{2\kappa}
=(2\kappa)^{-1}\e{L}(1;V)
L(1+\kappa,\chi_T;V\times V)W(2\kappa)
R^{2\kappa}+o(1),\eqno(4.30)
$$
by choosing
$P$ optimally and taking $\log R/\log T$ sufficiently large.
Then, on noting that the left side is larger
than $1$, we get $(4.25)$; and $(4.27)$ follows 
via $(1.2)$. In fact, it remains for us to show that 
$W(2\kappa)\asymp 1$; however, it should suffice to invoke
$W(2\kappa)=Y(1)$ and $(4.9)$.
\medskip
Further, we need to quote [6, Theorem 4][7, Vol.\ I,
(9.2.11)]:
\medskip
\noindent
{\bf Lemma 4.} {\it Let $v$ be a large positive parameter
and $\vartheta>0$ a constant. 
We put, with an integer $l\ge0$,
$$
\displaystyle{\Xi_d^{(l)}={1\over l!}(\vartheta\log v)^{-l}
\sum_{j=0}^l(-1)^{l-j}{l\choose j}\xi_d^{(j,l)},}\atop
\hbox{$\xi_d^{(j,l)}=\mu(d)
\Big({\log v^{1+j\vartheta}/d\Big)^l} $ 
for $d\le v^{1+j\vartheta}$, and  $= 0$ 
 for $d> v^{1+j\vartheta}$}.\eqno(4.31)
$$
Then we have that
$$
\hbox{$\Xi_d^{(l)} =\mu(d)$ for $d\le v$},\eqno(4.32)
$$
and 
$$
\sum_{n=1}^\infty d_l(n)
\Bigg(\sum_{d|n}\Xi_d^{(l)}\Bigg)^2n^{-\omega}\ll 1,
\eqno(4.33)
$$
whenever $ \omega\ge 1+1/\log v$.
Here $d_l(n)$ is the number of ways of representing 
$n$ as a product of $l$ positive integral factors, and the
implied constant depends on $l$ and $\vartheta$ 
at most.\/}
\bigskip
\noindent
{\bf 5. Linnik phenomenon extended.}
\par
\noindent
We proceed to the proof of
our main assertion which is given at the end of this section.
We begin with a conversion of [6, Lemma 5]
[7, Vol.\ I, (9.2.2)]
to our present situation: With
$\mu^2(r)=1$, we have, for $\Re s>1$,
$$
\eqalignno{
&\sum_{n=1}^\infty f(n)\tau_V(n)\Phi_r(n)
\Bigg(\sum_{d|n}\Xi^{(2)}_d\Bigg)n^{-s}
={L(s;V)L(s+\kappa,\chi_T;V)\over
L(2s+\kappa,\chi_T)}J_r(s),&(5.1)\cr
&J_r(s)={1\over K(r)}\sum_{d=1}^\infty 
\mu((d,r))\Xi^{(2)}_d
\prod_{p|d}\big(1-X_p(s)\big)
\prod_{\scr{p\nmid d}\atop\scr{p|r}}
\big(X_p(s)F_p-1\big),&(5.2)
}
$$
where $F_p$ is as in the previous section, and
$X_p(s)$ is the inverse of $(3.12)$ with $\chi=\chi_T$
and $a=-\kappa$. We set in Lemmas 2--4
$$
R=T^A,\; v=R^A,\; \vartheta=1/A,\eqno(5.3)
$$
with a sufficiently large constant $A>0$. 
Then, $J_r(s)$, $r\le R$, 
is regular and $J_r(s)\ll v^{1/2}$
for $\Re s\ge{3\over4}$.
With this, let us consider the expression
$$
I_r={1\over 2\pi i}\int_{2-i\infty}^{2+i\infty}
{L(\rho+w;V)L(\rho+w+\kappa,\chi_T;V)\over
L(2(\rho+w)+\kappa,\chi_T)}J_r(\rho+w)
\Gamma(w)Q^w dw.\eqno(5.4)
$$
Here $Q=v^{20}$, and
$$
L(\rho;V)=0,\> \rho=\beta+i\gamma;
\quad\txt{4\over5}\le\beta\le1,\;|\gamma|\le T.
\eqno(5.5)
$$
We have, by $(4.32)$ and $(5.1)$, 
$$
I_r=e^{-1/Q}+\sum_{v\le n}f(n)\tau_V(n)\Phi_r(n)
\Bigg(\sum_{d|n}\Xi^{(2)}_d\Bigg)n^{-\rho}e^{-n/Q}.
\eqno(5.6)
$$
We shift the line of integration in $(5.4)$ to 
$\Re w={3\over4}-\beta$ and get
$$
{1\over2}\le\Bigg|\sum_{v\le n\le  
Q^2}f(n)\tau_V(n)\Phi_r(n)
\Bigg(\sum_{d|n}\Xi_d^{(2)}\Bigg)n^{-\rho}
e^{-n/Q}\Bigg|^2,\eqno(5.7)
$$
since $I_r$ turns out to be negligibly small due to $(3.8)$ 
and the size of $J_r(s)$ 
mentioned above.
We multiply both sides by the factor 
$\mu^2(r)K(r)$ and sum over $r\le R$, getting
$$
\eqalignno{
&{G_1(R)\over\log T}\ll \sum_{\scr{N=2^\ell}
\atop\scr{v\le N\le Q^2}}
\sum_{r\le R}\mu^2(r)K(r)\cr
&\times\Bigg|\sum_{N\le n\le  
2N}f(n)\tau_V(n)\Phi_r(n)
\Bigg(\sum_{d|n}\Xi_d^{(2)}\Bigg)n^{-\rho}
e^{-n/Q}\Bigg|^2,&(5.8)
}
$$
where $\ell\in\B{N}$.
By virtue of Lemma 2 and $(4.27)$, we have
$$
{G_1(Z)\over\log T}\ll \e{L}(1;V)
L(1+\kappa,\chi_T;V\times V)Q^{4(1-\beta)}
\sum_{n=1}^\infty f(n)\Bigg(\sum_{d|n}\Xi_d^{(2)}
\Bigg)^2n^{-\omega_0},\eqno(5.9)
$$
with $\omega_0=1+(\log T)^{-1}$.
Then,  by $(4.33)$ with $l=2$ we find that
$$
{G_1(Z)\over \log T}\ll
 \e{L}(1;V)L(1+\kappa,\chi_T;V\times V)
Q^{4(1-\beta)}.\eqno(5.10)
$$
\par
Therefore, in view of $(4.25)$, we obtain 
\medskip
\noindent
{\bf Theorem.} {\it We assume the occurrence of 
the $T\!$-exceptional, among the zeros
of Dirichlet $L$-functions, in the sense $(1.1)$. 
Then, providing that the effectively computable 
absolute constant $a_0$ is 
adjusted appropriately, the inequality $(1.4)$
is uniformly satisfied by all
the zeros $\rho=\beta+i\gamma$, $|\gamma|\le T$,
of the Hecke $L$-function $L(s;V)$ associated with
any irreducible representation $V$ 
in the unitary principal series, with the spectral data
$|\nu_V|\le T$. \/}
\bigskip
\noindent
{\bf 6. Remarks.}
\par
\noindent
Some remarks are in order. The use of a quasi-character in 
the study of the zero density of Dirichlet $L$-functions was
initiated by Selberg [11], where actually the quasi-orthogonality
among Ramanujan sums was indicated. 
However, his character does not straightforwardly
generalise to $(4.12)$. To identify ours,
we need the observation by
the present author [5, p.\ 166] [6, p.\ 40] that Selberg's 
character originates
in fact in the $\Lambda^2$-sieve applied to the trivial
arithmetic function, i.e., the constant $1$. With
this, one may come to the idea $(4.3)$ and to 
the quasi-character $(4.12)$.
Nonetheless, the employment of $(4.3)$ is never trivial itself.
The insertion of the factor $\tau_V^2(n)$ is
deliberately made to exploit the Rankin convolution
of Maass cusp forms, which is
not necessary when dealing with holomorphic cusp forms
because of the validity of the Ramanujan; see [9]. We
stress that this idea generalises in conjunction
with the theory of the
Rankin convolution of symmetric power representations.
Such an instance is $(6.2)$ below.
Our present argument is of course applicable to the holomorphic
case as well, with minor changes; and our
theorem should hold 
regardless to which series of representations
a specific $V$ belongs. Also
our theorem holds for $L(s,\chi;V)$ with $\chi\bmod q$,
$\log q\ll\log T$, as well; the above restriction to the trivial
character is solely for the sake of simplicity. 
\par
One may deal with 
the zero-density of $L(s;V)$ by just extending
our large sieve inequality $(4.24)$ to an analogue of
[6, Lemma 4][7, Vol.\ I, (9.1.25) ]. 
It is possible also to hybridise $(4.24)$ 
by introducing the twist by primitive Dirichlet characters.
Then a large sieve zero density for sets of $L(s,\chi;V)$ 
with varying $\chi$
comes up that is an analogue of the assertion at the bottom
of [6, p.\ 196]. It should, however, be pointed out that
one may consider in general a single $V$ only, 
at least presently, even though the
end result can be made uniform in 
$V$ as our theorem manifests.
Nevertheless, if we restrict ourselves to the discrete series,
then the large sieve zero density involving both the sums
over $V$ and $\chi$ must be attainable; the reason for this
is that then the function $\Phi_r(n)$ 
can be made independent of $V$.
\par
It would be remiss not to mention that $L(s;V)$
for any irreducible $V$ 
does not have exceptional zeros in the sense that all the
zeros $\rho=\beta+i\gamma$ of $L(s;V)$ satisfy
$$
\beta<1-{c\over\log(|\nu_V|+|\gamma|)}
\eqno(6.1)
$$
with a computable absolute constant $c>0$. This important
result is due to Hoffstein--Ramakrishnan [3]; 
in fact, their assertion includes congruence subgroups as well.
Extensions of this to Rankin $L$-functions of two basic types
are achieved by Ramakrishnan--Wang [10]. These
are preceded by the pioneering work 
Goldfeld--Hoffstein--Lieman--Lockhart 
[2, Appendix] establishing the non-existence of
exceptional zeros for $L(s;\r{sym}^2V)$.
\par
As a matter of fact it appears to us to be more interesting
to apply the line of our argument to Rankin $L$-functions
than to Hecke $L$-functions. The aforementioned combination
of the $\Lambda^2$-sieve and the decomposition of the
Rankin convolution
$$
L(s;\r{sym}^2V\times\,\r{sym}^2V)=\zeta(s)
L(s;\r{sym}^2V)L(s;\r{sym}^4V),\eqno(6.2)
$$ 
which is an extension of $(3.14)$,
yields then the following analogues of 
Hoheisel's and Linnik's prime number theorems:
$$
\sum_{x\le p\le x+y}\tau_V^2(p)> C_1 y/\log x,
\quad |\nu_V|^{C_2}\le x, \quad x^{1-\theta}<y\le x;
\eqno(6.3)
$$
$$
\sum_{\scr{n\le x}\atop\scr{n\equiv a\bmod q}}
\tau_V^2(p)>C_3 x/(\varphi(q)\log x),\quad (a,q)=1,\quad
(q|\nu_V|)^{C_4}\le x,\eqno(6.4)
$$
where $\varphi$ is the Euler totient function. 
What is novel in these is that 
$\theta,\,C_j>0$ are all universal constants which
are effectively computable. 
Details will be developed in our forthcoming articles.
\vskip 1cm
\centerline{\bf References}
\medskip
\noindent
\item{[1]} S. Gelbart and H. Jacquet. A relation between
automorphic representations of $\r{GL}(2)$ and 
$\r{GL}(3)$. Ann.\ Sci.\ \'Ecole Normale Sup.\ $4^\r{e}$
s\'erie, {\bf 11} (1978), 471--552. 
\item{[2]} J. Hoffstein and P. Lockhart. Coefficients of Maass
forms and the Siegel zero. Ann.\ Math., {\bf 140} (1994), 
161--176; D. Goldfeld, J. Hoffstein and D. Lieman.
Appendix: An effective zero-free region. 
ibid, 177--181.
\item{[3]} J. Hoffstein and D. Ramakrishnan. Siegel zeros and
cusp forms. Intern.\ Math.\ Notices (1995), no.6, 279--308.
\item{[4]} H. Kim. Functionality for the exterior square of 
$\r{GL}_4$ and the symmetric fourth 
power of $\r{GL}_2$. J. Amer.\ Math.\
Soc., {\bf 16} (2003), 139--183.
\item{[5]} Y. Motohashi. Primes in arithmetic progressions.
Invent.\ math., {\bf 44} (1978), 163--178.
\item{[6]} ---. {\it Sieve Methods and 
Prime Number Theory\/}. 
Lect.\ Notes in Math.\ Phys., {\bf72},
Tata IFR and Springer-Verlag, Bombay 1983. 
http://www.math.tifr.res.in/$\sim$publ
/ln/tifr72.pdf
\item{[7]} ---. {\it Analytic Number Theory\/}.\ I. 
{\it Distribution of Prime Numbers\/}. Asakura Books, 
Tokyo 2009; II. {\it Zeta Analysis\/}. ibid, 2011:
ISBN 978-4-254-11821-6 and 978-4-254-11822-3,
respectively. (Japanese; translation into English 
being undertaken)
\item{[8]} ---. Elements of automorphic representations.
arXiv:1112.4226v1 [math.NT].
\item{[9]} ---. An extension of the Linnik phenomenon.
arXiv:1204.0149v1 [math.NT]; a revised version
to appear in Proc.\ Steklov Inst.\ Math.
\item{[10]} D. Ramakrishnan and S. Wang. 
On the exceptional
zeros of Rankin--Selberg $L$-functions. Compositio Math.,
{\bf 135} (2003), 211--244.
\item{[11]} A. Selberg. Remarks on sieves. In: 
Proc.\ 1972 Number Theory Conf., 
Boulder 1972, pp.\ 205--216. 
\item{[12]} G. Shimura. On the holomorphy of certain 
Dirichlet series. Proc.\ London Math.\ Soc., 
{\bf 31} (1975), 79--98.

\vskip 1cm
\noindent 
Honkomagome 5-67-1-901, Tokyo 113-0021, JAPAN
\hfill\def\ymzeta
{\font\brm=cmr17 at 30pt\font\sssrm=cmr5 at 4pt
\font\ssssrm=cmr5 at2.5pt
{{\brm O}\raise 9pt\hbox{\hskip -22pt
$\hfil\raise3pt\hbox{\ssssrm KH}\atop\hbox{
{\sssrm Y}$\zeta$\hskip-1pt{\sssrm M}}$}}}
\ymzeta

\bye